\documentstyle[12pt]{article}
\def\frak{\bf}
\begin{document}

\centerline{\bf Invariance of Plurigenera}

\bigskip
\centerline{Yum-Tong Siu\ %
\footnote{Partially supported by a grant from the National Science Foundation.}
}

\bigskip
\centerline{\it September 1997}

\bigskip
In this paper we give a proof of the following long conjectured
result on the invariance of the plurigenera.

\medskip
\noindent
{\it Main Theorem.} Let $\pi:X\rightarrow\Delta$ be a projective family
of compact complex manifolds parametrized by the open unit $1$-disk
$\Delta$.  
Assume that the family $\pi:X\rightarrow\Delta$ is of general type.
Then for every positive integer $m$
the plurigenus $\dim_{\bf C}\Gamma(X_t,mK_{X_t})$ is independent of
$t\in\Delta$, where $X_t=\pi^{-1}(t)$ and $K_{X_t}$ is the canonical
line bundle of $X_t$.

\medskip
\noindent
{\it Notations and Terminology.} The canonical line bundle of a complex
manifold $Y$ is denoted by $K_Y$.  The coordinate of the open unit $1$-disk
$\Delta$ is denoted by $t$.  Let $n$ be the complex dimension of
each $X_t$ for $t\in\Delta$.  In the assumption of the Main Theorem
the property of the family $\pi:X\rightarrow\Delta$ being of general type
means that for every $t\in\Delta$ there exist a positive integer $m_t$ and
a point $P_t\in X_t$ with the property that one can find elements
$s_0,s_1,\ldots,s_{n+1}\in\Gamma(X,m_tK_X)$ such that $s_0$ is nonzero 
at $P_t$ and ${s_1\over s_0},\ldots,{s_{n+1}\over s_0}$ form a 
local coordinate system of $X$ at $P_t$.  By the family 
$\pi:X\rightarrow\Delta$ being projective we mean that there exists 
a positive holomorphic line bundle on $X$.

\medskip
Let $K_{X,\pi}$ be the line bundle on $X$ 
whose restriction to $X_t$ is $K_{X_t}$ for each $t\in\Delta$.
Since the normal bundle of $X_t$ in $X$ is trivial, the two line
bundles $K_X$ and $K_{X,\pi}$ are naturally isomorphic. Under this
natural isomorphism a local section $s$ of $K_{X,\pi}$ corresponds to the 
local section $s\wedge\pi^*(dt)$ of $K_X$.  
Unless there is some risk of confusion, 
in this paper we will, without any further explicit mention, 
always identify $K_{X,\pi}$ with 
$K_X$ by this natural isomorphism.  Under this identification the
Main Theorem is equivalent to the statement that for every $t\in\Delta$
and every integer $m$ every element of $\Gamma(X_t,mK_{X_t})$ can 
be extended to an element of $\Gamma(X,mK_X)$.  

\medskip
The Hermitian metrics of holomorphic line bundles in this paper are allowed
to have singularities and may not be smooth.
For a Hermitian metric 
$h=e^{-\varphi}$ of a holomorphic line bundle $L$ over $X_0$ we denote
by ${\cal I}_h$ or by ${\cal I}_\varphi$ its multiplier ideal sheaf on
$X_0$ which by definition means the ideal sheaf on $X_0$
of all local holomorphic
function germs $f$ on $X_0$
such that $|f|^2e^{-\varphi}$ is integrable.  For the
proof of the Main Theorem, only multiplier ideal sheaves on $X_0$ are
considered and no multiplier ideal sheaves on $X$ are used.  
In the case of a Hermitian metric
$\tilde h=e^{-\tilde \varphi}$ of a holomorphic line bundle 
$\tilde L$ over $X$, for notational simplicity we simply use the notation
${\cal I}_{\tilde h}$ or ${\cal I}_{\tilde\varphi}$ 
to mean the multiplier ideal sheaf for the Hermitian metric
$\tilde h|X_0=e^{-\tilde \varphi}|X_0$ of the holomorphic line bundle
$\tilde L|X_0$ over $X_0$ and suppress the notation for
restriction to $X_0$.

\medskip
The stalk of a sheaf $\cal F$ at a point $P$ is denoted by ${\cal F}_P$.
The structure sheaf of a complex manifold $Y$ is denoted by ${\cal O}_Y$.
If $s$ is a global holomorphic section of a holomorphic line bundle
$L$ over a complex manifold $Y$ 
and if $\cal I$ is an ideal sheaf on $Y$, we say that
the germ of $s$ at a point $P$ belongs to ${\cal I}_P$ if 
the holomorphic function germ at $P$ which corresponds to the
germ of $s$ at $P$ with respect to some local trivialization of $L$
belongs to ${\cal I}_P$.
We say that $s$ is everywhere locally contained in $\cal I$ if 
at every point $P\in Y$ the germ of $s$ at $P$ belongs to ${\cal I}_P$,
or equivalently, $s$ belongs to the subset $\Gamma(Y,{\cal I}\otimes L)$
of $\Gamma(Y,L)$.  If $E$ is another holomorphic line bundle over $Y$
with a Hermitian metric $e^{-\chi}$, we say that $|s|^2e^{-\chi}$ is 
locally uniformly bounded (respectively locally integrable) on $Y$ 
if at every point $P\in Y$  
with respect to some local trivializations of $L$ and $E$ on some
open neighborhood $U$ of $P$ the function on $U$
corresponding to $|s|^2e^{-\chi}$ is uniformly bounded 
(respectively integrable) on $U$.  We say that $s$ is 
locally $L^2$ with respective to $e^{-\chi}$ on $Y$
if $|s|^2e^{-\chi}$ is locally integrable on $Y$.

\medskip
\noindent
{\it History and Sketch of the Proof of the Main Theorem.} 
Iitaka [I69-71] proved the special case of the invariance of the
plurigenera in a family of surfaces.  His method works only for
surfaces because it uses much of the information 
from the classification of surfaces.
Levine [L83,L85] proved that for every positive integer $m$
every element of
$\Gamma(X_0,mK_{X_0})$ can be extended to the fiber of $X$ over the
double point of $\Delta$ at $t=0$.  So far there is no way to
continue the process to get an extension to the fiber of $X$
over a point of $\Delta$ at $t=0$ of any finite order.  
Nakayama [Nak86] pointed out that if the relative case of the
minimal model program can be carried out for a certain dimension, 
the conjecture of the invariance of the plurigenera for that
dimension would follow directly from it.  Thus the invariance of the 
plurigenera for threefolds is a consequence of the completion of 
the relative case of the minimal model program for the case of threefolds 
by Kollar and Mori [KM92].

\medskip
For the proof of the Main Theorem here we use a strategy
completely different from those used by the others in the past.  We 
now sketch our strategy and leave out the less essential technical
details.  There are some unavoidable technical inaccuracies in the sketch
due to the suppression of precise details.
There are three ingredients in our proof: Nadel's multiplier
ideal sheaves [Nad89], Skoda's result on the generation of ideals with
$L^2$ estimates with respect to a plurisubharmonic weight [Sk72], and the 
extension theorem of Ohsawa-Takegoshi-Manivel for holomorphic
top-degree forms which are $L^2$ with respect to a plurisubharmonic
weight [OT87,M93].  The extension theorem of
Ohsawa-Takegoshi-Manivel is for the setting of a Stein domain or
manifold and a global plurisubharmonic function as weight.  Here
we adapt it to the case of a projective family of compact complex manifolds
and a Hermitian metric of a line bundle with nonnegative curvature current.  
The adaptation is done by restricting to a Stein Zariski open subset 
on which the line
bundle is globally trivial, because $L^2$ bounds automatically extend
the domain of definition from the Zariski open subset to 
the family of compact manifolds.

\medskip
We take the $m^{\hbox{th}}$ roots of basis elements
of $\Gamma(X_0,mK_{X_0})$ for every positive integer $m$ to use them
in an infinite series to construct
a Hermitian metric $e^{-\varphi}$ for $K_{X_0}$.
We also take
the $m^{\hbox{th}}$ roots of basis elements
of $\Gamma(X,mK_X)|X_0$ for every positive integer $m$ to use them
in an infinite series to construct
a Hermitian metric $e^{-\tilde\varphi}$ for $K_{X_0}$.  If for an
infinite number of integers $\ell_\nu$ the singularity of
$e^{-\ell_\nu\tilde\varphi}$ are only worse than that of
$e^{-\ell_\nu\varphi}$ by
some fixed amount independent of $\nu$, then the
extension theorem of Ohsawa-Takegoshi-Manivel can be applied to
yield the Main Theorem.

\medskip
If the contrary holds, then for some appropriate positive integer $m$
we can construct some Hermitian metric for $mK_{X_0}$
whose singularity is suitably chosen to be between those of  
$e^{-m\tilde\varphi}$ and $e^{-m\varphi}$ so that, by Skoda's result,
we can use this Hermitian metric to produce an
element $s$ of $\Gamma(X_0,mK_{X_0})$ which is 
$L^2$ with respect to $e^{-(m-1)\tilde\varphi}$ but not locally $L^2$ with
respect to $e^{-m\tilde\varphi}$ everywhere on $X_0$.  On the other hand,
by the extension theorem of Ohsawa-Takegoshi-Manivel we can regard $s$
as an $(m-1)K_{X_0}$-valued top-degree form on $X_0$ which is
$L^2$ with respect to $e^{-(m-1)\tilde\varphi}$ and 
can therefore extend it
to an element of $\Gamma(X,mK_X)$.  The definition of
$e^{-\tilde\varphi}$ implies that $|s|^2\, e^{-m\tilde\varphi}$ is locally
uniformly bounded on $X_0$
and consequently $s$ is $L^2$ with respect to 
$e^{-m\tilde\varphi}$ everywhere on $X_0$, which is a contradiction.

\medskip
One of the technical details is that the Hermitian
metric $e^{-m\tilde\varphi}$
has to be slightly modified to make sure that 
its curvature current dominates a
smooth positive $(1,1)$-form in order to apply 
the extension theorem of Ohsawa-Takegoshi-Manivel.  For that modification
the Kodaira technique of 
writing some high multiple of a big line bundle as an effective divisor
plus an ample line bundle is used.

\medskip
\noindent
{\it Lemma 1.} Let $L$ be a holomorphic line bundle over an
$n$-dimensional compact complex manifold $Y$ with a
Hermitian metric which is locally of the form $e^{-\xi}$ with
$\xi$ plurisubharmonic.  Let ${\cal I}_\xi$ be the multiplier
ideal sheaf of the Hermitian metric $e^{-\xi}$.  Let $E$ be an ample holomorphic
line bundle over $Y$ such that for every point $P$ of $Y$ there are a
finite number of elements of $\Gamma(Y,E)$ which 
all vanish to order at least $n+1$ at $P$
and which do not simultaneously vanish outside $P$.
Then $\Gamma(Y,{\cal I}_\xi\otimes(L+E+K_Y))$ generates
${\cal I}_\xi\otimes(L+E+K_Y)$ at every point of $Y$.

\medskip
\noindent
{\it Proof.} The key ingredient is the following result of Skoda 
[Sk72, Th.1, pp.555-556].

\medskip
Let $\Omega$ be a pseudoconvex 
domain in ${\bf C}^n$ and $\psi$ be a plurisubharmonic function on 
$\Omega$.  Let
$g_1,\ldots,g_p$ be holomorphic functions on $\Omega$.  Let
$\alpha>1$ and $q=\inf(n,p-1)$.  Then for every holomorphic function $f$ on
$\Omega$ such that
$$
\int_\Omega|f|^2|g|^{-2\alpha q-2}e^{-\psi}d\lambda<\infty,
$$
there exist holomorphic functions $h_1,\ldots,h_p$ on $\Omega$ such that
$$
f=\sum_{j=1}^pg_jh_j
$$
and
$$
\int_\Omega|h|^2|g|^{-2\alpha q}e^{-\psi}d\lambda
\leq{\alpha\over{\alpha-1}}\int_\Omega|f|^2|g|^{-2\alpha q-2}e^{-\psi}d\lambda,
$$
where
$$
|g|=\left(\sum_{j=1}^p|g_j|^2\right)^{1\over 2},\qquad
|h|=\left(\sum_{j=1}^p|h_j|^2\right)^{1\over 2},
$$
and $d\lambda$ is the Euclidean volume element of ${\bf C}^n$.

\medskip
Fix arbitrarily $P_0\in Y$.
Take an arbitrary element $s$ of $({\cal I}_\xi)_{P_0}$.
Let $z=(z_1,\ldots,z_n)$ be a local coordinate system on some open
neighborhood $U$ of $P_0$ with $z(P_0)=0$ such that $L|U$ is trivial.
Let $\rho$ be a cut-off function centered at $P_0$ so that $\rho$ is a
smooth nonnegative-valued function with compact support in $U$ which
is identically $1$ on some Stein open neighborhood $\Omega$ of $P_0$.  Choose
$u_1,\ldots,u_N\in\Gamma(Y,E)$
whose common zero-set consists of the single point $P_0$ and which
all vanish to order at least $n+1$ at $P_0$.
Let $h_E$ be a smooth Hermitian metric of $E$ whose curvature form is
strictly positive at every point of $Y$.  Let $0<\eta<{1\over n+1}$.
By the standard techniques
of $L^2$ estimates of $\bar\partial$, we can solve the equation
$$
\bar\partial\sigma=\rho\bar\partial s
$$
for a smooth section $\sigma$ of $L+E+K_Y$ over $Y$
which is $L^2$ with respect to the Hermitian metric
$$
e^{-\xi}\left(h_E\right)^\eta\over
\left(\sum_{j=1}^N|u_j|^2\right)^{1-\eta}
$$
of $L+E$.
Then $\rho s-\sigma$ is an element of 
$\Gamma(Y,{\cal I}_\xi\otimes(L+E+K_Y))$.
Since $\rho\bar\partial s$ is identically zero on $\Omega$, it follows that
$\sigma$ is holomorphic on $\Omega$.
We now apply Skoda's result to the case $g_j=z_j$ ($1\leq j\leq n$) with
$q=n-1$ and $\alpha={(1-\eta)(n+1)-1\over n-1}>1$
and $\psi=\xi$.  (For the case $n=1$ we simply choose $\alpha$ be any
number greater than $1$, because in that case $\alpha q$ is always zero.)
Let $|z|=\left(\sum_{j=1}^n|z_j|^2\right)^{1\over 2}$.
Since $u_1,\ldots,u_N$ all vanish to order at least $n+1$ at $P_0$, it 
follows that
$$
\int_\Omega|\sigma|^2e^{-\xi}|z|^{-2\alpha q-2}=
\int_\Omega|\sigma|^2e^{-\xi}|z|^{-2(1-\eta)(n+1)}<\infty.
$$
By Skoda's reuslt
$$
\sigma=\sum_{j=1}^n\tau_j z_j
$$
locally at $P_0$ for some 
$\tau_1,\ldots,\tau_n\in\left({\cal I}_\xi\right)_{P_0}$.

\medskip
Let $J$ be the ideal at $P_0$ generated by elements of
$\Gamma(Y,{\cal I}_\xi\otimes(L+E+K_Y))$ over $\left({\cal O}_Y\right)_{P_0}$.  
It follows from
$$
\rho s-\sigma\in\Gamma(Y,{\cal I}_\xi\otimes(L+E+K_Y))
$$
that
$$
s\in J+{\frak m}_{P_0}\left({\cal I}_\xi\right)_{P_0},
$$
where ${\frak m}_{P_0}$ is the maximum ideal of $Y$ at $P_0$.
Since $s$ is an arbitrary element of $\left({\cal I}_\xi\right)_{P_0}$, it 
follows that
$$
\left({\cal I}_\xi\right)_{P_0}\subset 
J+{\frak m}_{P_0}\left({\cal I}_\xi\right)_{P_0}.
$$
Clearly we have $J\subset\left({\cal I}_\xi\right)_{P_0}$.  Thus
$$
\left({\cal I}_\xi\right)_{P_0}/J\subset{\frak m}_{P_0}
\left(\left({\cal I}_\xi\right)_{P_0}/J\right).
$$
By Nakayama's lemma,
$$
\left({\cal I}_\xi\right)_{P_0}/J=0
$$
and $J=\left({\cal I}_\varphi\right)_{P_0}$. Q.E.D.

\bigskip
For every positive integer $m$ we choose a basis
$$
s^{(m)}_1,\cdots,s^{(m)}_{q_m}\in\Gamma\left(X_0,mK_{X_0}\right)
$$
and we choose
$$
\tilde s^{(m)}_1,\cdots,\tilde s^{(m)}_{\tilde q_m}
\in\Gamma\left(X,mK_X\right)
$$
so that
$$
\tilde s^{(m)}_1\big|{X_0},\cdots,\tilde s^{(m)}_{\tilde q_m}\big|{X_0}
$$
is a basis of $\Gamma\left(X,mK_X\right)\big|_{X_0}$
and $\tilde s^{(m)}_\nu=
s^{(m)}_\nu$ for $1\leq\nu\leq\tilde q_m$.
We can choose a sequence of positive numbers $\theta_m$ so that
$$
\sum_{m=1}^\infty\theta_m\left(\sum_{\nu=1}^{q_m}\left|s^{(m)}_\nu\right|
^{2\over m}\right)
$$
converges uniformly on compact subsets of $X_0$ to a Hermitian
metric of $-K_{X_0}$ and
$$
\sum_{m=1}^\infty\theta_m\left(\sum_{\nu=1}^{\tilde q_m}\left|
\tilde s^{(m)}_\nu\right|
^{2\over m}\right)
$$
converges uniformly on compact subsets of $X$ to a Hermitian metric
of $-K_X$.
Locally on $X_0$ we define
$$
\varphi=\log\sum_{m=1}^\infty\theta_m\left(\sum_{\nu=1}^{q_m}\left|s^{(m)}_\nu\right|
^{2\over m}\right)
$$
so that $e^{-\varphi}$ is a Hermitian metric of $K_{X_0}$.
Locally on $X$ we define
$$
\tilde\varphi=\log\sum_{m=1}^\infty\theta_m
\left(\sum_{\nu=1}^{\tilde q_m}\left|
\tilde s^{(m)}_\nu\right|
^{2\over m}\right)
$$
so that $e^{-\tilde\varphi}$ is a Hermitian metric of $K_X$.

\medskip
Since the family $\pi:X\rightarrow\Delta$ is of general type, 
we can choose an integer 
$m_0\geq 2$ such that
$m_0K_X=D+F$, where $D$ is an effective divisor on $X$ not containing
$X_0$ and $F$ is such a high multiple of a positive line bundle on $X$ 
that

\medskip
\noindent
(i) for every point $P\in X_0$ there exist a finite number of
elements of $\Gamma(X,F-2K_X)|X_0$ whose common zero-set consists
only of the single point $P$ 
and which all vanish to order at least $n+1$ at $P$ and

\medskip
\noindent
(ii) a basis of $\Gamma(X,F)|X_0$ embeds $X_0$ as a complex
submanifold of some complex projective space.

\medskip
\noindent
Let $s_D$ be the canonical
section of the holomorphic line bundle $D$ so that the divisor of $s_D$
is $D$.  Let
$$
u_1,\ldots,u_N\in\Gamma(X,F)
$$
such that
$$
u_1\bigm|X_0,\ldots,u_N\bigm|X_0
$$
form a basis of $\Gamma(X,F)|X_0$.
From $s_Du_j\in\Gamma(X,m_0K_X)$ ($1\leq j\leq N$) and 
the non simultaneous vanishing of $u_1,\ldots,u_N$ at any point of $X_0$
it follows from the definition of $\tilde\varphi$
that $|s_D|^2e^{-m_0\tilde\varphi}|_{X_0}$ is locally uniformly bounded 
on $X_0$.  Let
$$
h_F={1\over\sum_{j=1}^N\left|u_j\right|^2}
$$
and we introduce the Hermitian metric
$$
e^{-\psi}=\left({h_F\over|s_D|^2}\right)^{1\over m_0}
$$
for the line bundle $K_X$.  Choose $0<\epsilon<1$ such that
$e^{-\epsilon\psi}|X_0$ is locally integrable on $X_0$.
For any positive integer $\ell$ we introduce the Hermitian metric
$$
h_\ell=e^{-(\ell-\epsilon)\varphi-(m_0+\epsilon)\psi}
$$
for the line bundle $(\ell+m_0)K_{X_0}$.

\medskip
As stated at the beginning of the paper, in the statement of Lemma 2 
below and for the rest of the paper the notation
${\cal I}_{(\ell_\nu+m_0-1-\epsilon)\tilde\varphi+\epsilon\psi}$
denotes an ideal sheaf on $X_0$ and not an ideal sheaf on $X$ and it is
the multiplier ideal sheaf for the metric 
$e^{-(\ell_\nu+m_0-1-\epsilon)\tilde\varphi-\epsilon\psi}|X_0$
of the holomorphic line bundle $(\ell_\nu+m_0-1)K_{X_0}$.

\medskip
\noindent
{\it Lemma 2.} Let $\ell_0$ be a positive integer.
Suppose there exists a sequence of positive integers
$\ell_\nu\nearrow\infty$ ($1\leq\nu<\infty$)
such that
$$
{\cal I}_{h_{\ell_\nu}}\subset{\cal I}
_{(\ell_\nu+m_0-1-\epsilon)\tilde\varphi+\epsilon\psi}.
$$
Then any element $s$ of $\Gamma(X_0,\ell_0K_{X_0})$ is everywhere locally
contained in ${\cal I}_{\ell_0\tilde\varphi}$.

\medskip
\noindent
{\it Proof}.  Without loss of generality we can assume after reindexing
the sequence $\{\ell_\nu\}_{1\leq\nu<\infty}$
that $\ell_\nu>2\ell_0$
for $1\leq\nu<\infty$.  Take an arbitrary $P_0\in X_0$.  Let
$\ell_\nu=q_\nu\ell_0+r_\nu$ with $0\leq r_\nu<\ell_0$ ($1\leq\nu<\infty$).
Take a non-identically-zero
$\sigma\in\left({\cal O}_X\right)_{P_0}$ such that $|\sigma|^2e^{-\ell_0\tilde\varphi}$
is bounded in the supremum norm on some open neighborhood $U$ of $P_0$
in $X_0$ (for example, we can take $\sigma$ to be
the germ at $P_0$ of some nonzero element of $\Gamma(X,\ell K_X)|X_0$ for 
some $\ell\geq\ell_0$).  Since 
$\left|s^{q_\nu}\right|^2e^{-q_\nu\ell_0\varphi}$
is uniformly bounded on $X_0$ from the definition of $\varphi$, it follows from
$0\leq r_\nu<\ell_0$ and $\tilde\varphi\leq\varphi$ and the
integrability of $e^{-\epsilon\psi}$ that the germ of $s^{q_\nu}\sigma s_D$
at $P_0$ belongs to $\left({\cal I}_{h_{\ell_\nu}}\right)_{P_0}$.
By the assumption of the Lemma, the germ of
$s^{q_\nu}\sigma s_D$ at $P_0$ belongs to
$\left({\cal I}
_{(\ell_\nu+m_0-1-\epsilon)\tilde\varphi+\epsilon\psi}\right)_{P_0}$.
There exists some relatively compact open neighborhood 
$W$ of $P_0$ in $U$ with $K_{X_0}|W$ trivial such that
$$
\int_W\left|s^{q_\nu}\sigma s_D\right|^2
e^{-(\ell_\nu+m_0-1-\epsilon)\tilde\varphi-\epsilon\psi}<\infty.
$$

\medskip
Let ${1\over q_\nu}+{1\over q^\prime_\nu}=1$.  Then
$q^\prime_\nu={q_\nu\over q_\nu-1}$ and 
${q^\prime_\nu\over q_\nu}={1\over q_\nu-1}$ and we have by H\"older's
inequality
$$
\displaylines{
\int_W|s|^2e^{-\ell_0\tilde\varphi}=
\int_W\left|s\sigma^{1\over q_\nu}s_D^{1\over q_\nu}
\sigma^{-1\over q_\nu}s_D^{-1\over q_\nu}\right|^2e^{-\ell_0\tilde\varphi}
e^{-\epsilon\psi\over q_\nu}e^{\epsilon\psi\over q_\nu}\cr
\leq
\left(\int_W\left|s^{q_\nu}\sigma s_D\right|^2
e^{-q_\nu\ell_0\tilde\varphi-\epsilon\psi}\right)^{1\over q_\nu}
\left(\int_W\left|\sigma^{-q^\prime_\nu\over q_\nu}
s_D^{-q^\prime_\nu\over q_\nu}\right|^2
e^{\epsilon\psi q^\prime_\nu\over q_\nu}\right)^{1\over q^\prime_\nu}\cr
=
\left(\int_W\left|s^{q_\nu}\sigma s_D\right|^2e^{-q_\nu\ell_0\varphi
-\epsilon\psi}\right)^{1\over q_\nu}
\left(\int_W\left|\sigma^{-1\over q_\nu-1}
s_D^{-1\over q_\nu-1}\right|^2e^{\epsilon\psi\over q_\nu-1}
\right)^{q_\nu\over q_\nu-1}\cr
\leq
C\left(\int_W\left|s^{q_\nu}\sigma s_D\right|^2
e^{-(\ell_\nu+m_0-1-\epsilon)\tilde\varphi-\epsilon\psi}\right)^{1\over q_\nu}
\left(\int_W\left|\sigma^{-1\over q_\nu-1}
s_D^{-1\over q_\nu-1}\right|^2e^{\epsilon\psi\over q_\nu-1}\right)^{q_\nu\over q_\nu-1},\cr
}
$$
where $C$ is the supremum of
$e^{(m_0+r_\nu-1)\tilde\varphi\over q_\nu}$ on $W$.
For $q_\nu$ sufficiently large,
$$
\int_W\left|\sigma^{-1\over q_\nu-1}
s_D^{-1\over q_\nu-1}\right|^2e^{\epsilon\psi\over q_\nu-1}<\infty.
$$
Hence
$$
\int_W|s|^2e^{-\ell_0\tilde\varphi}<\infty
$$
and the germ of $s$ at $P_0$ belongs 
to $\left({\cal I}_{\ell_0\tilde\varphi}\right)_{P_0}$.
Q.E.D.

\medskip
For the next step in our proof of the Main Theorem we need the following
extension statement which is an adaptation of the extension theorem of 
Ohsawa-Takegoshi-Manivel.

\medskip
\noindent
{\it Proposition 3}.  Let $\gamma:Y\rightarrow\Delta$ be a projective family
of compact complex manifolds parametrized by the open unit $1$-disk
$\Delta$. Let $Y_0=\gamma^{-1}(0)$ and let
$n$ be the complex dimension of $Y_0$. Let
$L$ be a holomorphic line bundle with a Hermitian
metric which locally is represented by $e^{-\chi}$ such that
$\sqrt{-1}\partial\bar\partial\chi\geq \omega$ in the sense
of currents for some smooth
positive $(1,1)$-form $\omega$ on $Y$.
Let $0<r<1$ and $\Delta_r=\{\,t\in\Delta\bigm|\,|t|<r\,\}$.  Then there 
exists a positive constant $A_r$ with the following property.  For any
holomorphic $L$-valued $n$-form $f$ on $Y_0$ with
$$
\int_{Y_0}|f|^2e^{-\chi}<\infty,
$$
there exists a holomorphic $L$-valued $(n+1)$-form $\tilde f$ on 
$\gamma^{-1}(\Delta_r)$
such that $\tilde f|_{Y_0}=f\wedge\gamma^*(dt)$ at points of $Y_0$ and
$$
\int_Y|\tilde f|^2e^{-\chi}\leq A_r
\int_{Y_0}|f|^2e^{-\chi}.
$$
Here no metrics of the tangent bundles of $Y_0$ and $Y$ are needed
to define the integrals of the absolute-value squares of top-degree 
holomorphic forms
$f$ and $\tilde f$ respectively on $Y_0$ and $Y$.

\medskip
\noindent
{\it Proof.} The proof can be easily adapted in the following way from 
the techniques given in [Si96] for the alternative
proof there of the theorem of Ohsawa-Takegoshi. 
(Proofs can also be obtained by
modifying those in [OT83, M93].) Let 
$v$ be a meromorphic section of $L$ over $Y$ 
so that neither the pole-set nor the zero-set of $v$
contains $Y_0$.  Choose a complex hypersurface $Z$ in $Y$
containing the zero-set and the pole-set of $v$ such that $Z$ does not
contain $Y_0$ and $Y-Z$ is Stein.  For every positive integer $\nu$
let $\Omega_\nu$ be a relatively compact Stein open subset of $X-Z$ with
smooth strictly pseudoconvex boundary such that 
$\cup_{\nu=1}^\infty\Omega_\nu=X-Z$ and 
$\Omega_{\nu}$ is relatively compact in $\Omega_{\nu+1}$.
On $X-Z$ under the isomorphism defined by 
division by $v$ the line bundle $L|X-Z$ is globally trivial.  We let
$\tilde\chi$ be
the plurisubharmonic function $-\log\left(|v|^2e^{-\chi}\right)$
on $X-Z$.  

\medskip
We now apply the techniques in [Si96]
of the alternative proof of the theorem of
Ohsawa-Takegoshi to extend, after multiplication by $\gamma^*(dt)$,
the top-degree holomorphic form ${f\over v}$ on
$\Omega_\nu\cap Y_0$ which is $L^2$ 
on $\Omega_\nu\cap Y_0$ with respect to $e^{-\tilde\chi}$ 
to a top-degree holomorphic form $G_\nu$ 
on $\gamma^{-1}(\Delta_r)\cap\Omega_\nu$
whose $L^2$ norm on $\gamma^{-1}(\Delta_r)\cap\Omega_\nu$
with respect to $e^{-\tilde\chi}$ is bounded by a finite constant independent
of $\nu$.  When we apply the techniques of the alternative
proof of the theorem of 
Ohsawa-Takegoshi, we have to use holomorphic tangent vector fields of
the Stein manifold $\Omega_{\nu+1}$ to get a 
sequence of smooth plurisubharmonic
functions on $\Omega_\nu$ which approach the plurisubharmonic
function $\tilde\chi$ on $\Omega_\nu$.  The extension $\tilde f$ 
is obtained as
the limit of $G_\nu v$ as $\nu$ goes to infinity.  The 
smooth positive $(1,1)$-form $\omega$ in the assumption is needed
for the $\nu$-independent {\it a priori} estimates 
for the solution of the modified $\bar\partial$ equation on 
$\gamma^{-1}(\Delta_r)\cap\Omega_\nu$ in the techniques of the 
alternative proof of the
theorem of Ohsawa-Takegoshi.  Q.E.D.

\medskip
\noindent
{\it Lemma 4}.  If $m$ is an integer $\geq 2$ and if
$s\in\Gamma(X_0,mK_{X_0})$ is everywhere locally contained in 
${\cal I}_{(m-1-\epsilon)\tilde\varphi+\epsilon\psi}$,
then $s$ can be extened to an element of $\Gamma(X,mK_X)$.

\medskip
\noindent
{\it Proof.}
We apply Proposition 3 to the case $L=(m-1)K_X$, 
$\chi=(m-1-\epsilon)\tilde\varphi+\epsilon\psi$, and $f=s$ to extend
$s$ to an $(m-1)K_X$-valued holomorphic $(n+1)$-form on $\pi^{-1}(\Delta_r)$
for some $0<r<1$, where $\Delta_r$ is the open 1-disk of radius $r$
centered at the origin. Then we use the standard theory of coherent sheaves
and Stein spaces to get the extension from $\pi^{-1}(\Delta_r)$ to
all of $X$.  Q.E.D.

\medskip
\noindent
{\it Lemma 5}. If $m$ is an integer $\geq 2$ and if 
$s\in\Gamma(X_0,mK_{X_0})$ is everywhere locally contained in 
${\cal I}_{m\tilde\varphi}$, then
$s$ can be extended to an element
of $\Gamma(X,mK_X)$.

\medskip
\noindent
{\it Proof.} Since $s$ is everywhere locally contianed in
${\cal I}_{m\tilde\varphi}$, we can cover $X_0$ by a finite number of
open subsets $U_j$ ($1\leq j\leq k$) such that
$K_{X_0}|U_j$ is trivial on $U_j$ and
$$
\int_{U_j}|s|^2e^{-m\tilde\varphi}<\infty
$$
for $1\leq j\leq k$.
Take $0<\eta<1$ and consider the Hermitian metric 
$e^{-(m-1-\eta)\tilde\varphi-\eta\psi}$
for $(m-1)K_X$.  We apply H\"older's inequality with $p={m\over m-1-\eta}$ 
and $p^\prime={p\over p-1}={m\over 1+\eta}$ to get
$$
\int_{U_j}|s|^2
e^{-(m-1-\eta)\tilde\varphi-\eta\psi}
\leq
\left(\int_{U_J}|s|^2
e^{-m\tilde\varphi}\right)^{m-1-\eta\over m}
\left(\int_{U_j}|s|^2e^{-m\eta\psi\over 1+\eta}
\right)^{1+\eta\over m}.
$$
When $\eta$ is sufficiently small,
$$
e^{-m\eta\psi\over 1+\eta}=\left({h_F\over|s_D|^2}\right)
^{-m\eta\psi\over (1+\eta)m_0}
$$
is locally integrable at every point of $X_0$.  Hence $s$ is $L^2$
as an $(m-1)K_{X_0}$-valued $n$-form on
$X_0$ with respect to the Hermitian metric
$$
e^{-(m-1-\eta)\tilde\varphi-\eta\psi}|X_0
$$
whose curvature current, because of the factor $h_F$, is bounded from below
by a smooth positive $(1,1)$-form on $X_0$.  Now the Lemma follows from 
Proposition 3 (or Lemma 4).  Q.E.D.

\medskip
\noindent
{\it Final Step of the Proof of the Main Theorem.}  
From the definition of $h_\ell$ for $\ell=1$ we have
$$
{\cal I}_{h_1}={\cal I}_{(1-\epsilon)\varphi+(m_0+\epsilon)\psi}
\subset{\cal I}_{(m_0-\epsilon)\tilde\varphi+\epsilon\psi},\leqno{(1)}
$$
because $|s_D|^2e^{-m_0\tilde\varphi}|_{X_0}$ is locally uniformly bounded
on $X_0$.
Fix an arbitrary positive integer $\ell_0$.
To prove the Main Theorem, it 
suffices to show that every element of $\Gamma(X_0,\ell_0 K_{X_0})$ can be
extended to an element of $\Gamma(X,\ell_0 K_X)$.
Suppose the contrary and we are going to derive a contradiction.
By Lemma 2 and Lemma 5 we can assume that
there exists a positive integer $\ell_\#$ 
such that for $\ell\geq\ell_\#$ we have
$$
{\cal I}_{h_\ell}\not\subset{\cal I}
_{(\ell+m_0-1-\epsilon)\tilde\varphi+\epsilon\psi}.
$$
By (1) we know that there is a smallest positive
integer $\ell_*$ (which must be at least $2$) such that
$$
{\cal I}_{h_{\ell_*}}\not\subset{\cal I}
_{(\ell_*+m_0-1-\epsilon)\tilde\varphi+\epsilon\psi}.\leqno{(2)}
$$
Then
$$
{\cal I}_{h_{\ell_*-1}}\subset{\cal I}
_{(\ell_*-1+m_0-1-\epsilon\psi)\tilde\varphi+\epsilon\psi}.\leqno{(3)}
$$
From the choice of $F$ we know that
the line bundle $(F|X_0)-K_{X_0}$ over $X_0$ is ample.
We now apply Lemma 1 to the case
$E=(F|X_0)-2K_{X_0}$ and $L=\ell_*K_{X_0}+D$ with the Hermitian metric
$$
e^{-\xi}={e^{-(\ell_*-\epsilon)\varphi-\epsilon\psi}\over|s_D|^2}.
$$
Each of the two Hermitian 
metrics $h_{\ell_*}$ and $e^{-\xi}$ is locally bounded on $X_0$ by
a positive constant times the other.  Hence
the two ideal sheaves
${\cal I}_\xi$ and ${\cal I}_{h_{\ell_*}}$ coincide everywhere
on $X_0$.
By Lemma 1 it follows from $(\ell_*+m_0-1)K_{X_0}=L+E+K_{X_0}$ that 
$\Gamma(X_0,{\cal I}_{h_{\ell_*}}\otimes(\ell_*+m_0-1)K_{X_0})$ locally
generates ${\cal I}_{h_{\ell_*}}$ on $X_0$.  From (2) it follows that
there exists
$s\in\Gamma(X_0,{\cal I}_{h_{\ell_*}}\otimes(\ell_*+m_0-1)K_{X_0})$
such that 
$$ s\ \hbox{ is not everywhere locally contained in } 
\ {\cal I}_{(\ell_*+m_0-1-\epsilon)\tilde\varphi+\epsilon\psi}.\leqno{(4)}
$$
From (3) and ${\cal I}_{h_{\ell_*}}\subset {\cal I}_{h_{\ell_*-1}}$ 
it follows that every element of 
$\Gamma(X_0,{\cal I}_{h_{\ell_*}}\otimes(\ell_*+m_0-1)K_{X_0})$ 
is locally contained in
${\cal I}_{(\ell_*-1+m_0-1-\epsilon)\tilde\varphi+\epsilon\psi}$
at every point of $X_0$.  
From Lemma 4 it follows that
$s$ can be extended to an element $\tilde s$ of 
$\Gamma(X,(\ell_*+m_0-1)K_X)$.  Since from the definition
of $\tilde\varphi$ we know that 
$|\tilde s|^2e^{-(\ell_*+m_0-1)\tilde\varphi}$
is uniformly bounded on $X_0$,
it follows from the integrability of $e^{-\epsilon\psi}$
that $s$ is everywhere locally contained in 
${\cal I}_{(\ell_*+m_0-1-\epsilon)\tilde\varphi+\epsilon\psi}$,
which contradicts (4).
This concludes the proof of the Main Theorem.  Q.E.D.

\bigskip
\noindent
{\it References.}

\medskip
\noindent
[I69-71] S. Iitaka, Deformations of compact complex surfaces I, II, and III.
In: {\it Global Analysis}, papers in honor of K. Kodaira, Princeton University
Press 1969, pp.267-272; J. Math. Soc. Japan 22 (1970), 247-261; {\it ibid}
23 (1971), 692-705.

\noindent
[KM92] Kollar, J., Mori, S.: Classification of three dimensional flips,
Journal of Amer. Math. Soc. 5 (1992), 533-702.

\noindent
[L83] Levine, M.: Pluri-canonical divisors on K\"ahler manifolds, Invent. math.
74 (1983), 293-903.

\noindent
[L85] Levine, M.: Pluri-canonical divisors on K\"ahler manifolds II, Duke 
Math. J. 52 (1985), 61-65.

\noindent
[M93] Manivel, L.: Un th\'or\`eme de prolongement $L^2$ de sections 
holomorphes d'un fibr\'e hermitien. Math. Zeitschr. 212, 107-122 (1993).

\noindent
[Nad89] Nadel, A.: Multiplier ideal sheaves and the existence of K\"ahler-Einstein
metrics of positive scalar curvature, Proc. Natl. Acad. Sci. USA, 86, 
7299-7300 (1989), and Ann. of Math., 132, 549-596 (1989).

\noindent
[Nak86] Nakayama, N.: Invariance of plurigenera of algebraic varieties 
under minimal model conjecture, Topology 25 (1986), 237-251. 

\noindent
[OT87] Ohsawa, T., Takegoshi, K.: On the extension of  $L^2$ holomorphic 
functions, Math. Z., 195, 197-204 (1987).

\noindent
[Si96] Siu, Y.-T., The Fujita conjecture and the extension theorem of 
Ohsawa-Takegoshi, in {\it Geometric Complex Analysis} 
ed. Junjiro Noguchi {\it et al},
World Scientific: Singapore, New Jersey, London, Hong Kong 1996, pp. 577-592.

\noindent
[Sk72] H  Skoda, Application des techniques $L^2$ \`a la th\'eorie des ideaux d'un 
alg\`ebre de fonctions holomorphes avec poids, Ann. Sci. Ec. Norm. Sup. 5 
(1972), 548-580. 

\medskip
\noindent
Author's address:

\medskip
\noindent
Yum-Tong Siu, Department of Mathematics, Harvard University, Cambridge, 
MA 02193
(e-mail: siu@math.harvard.edu).

\end{document}